# Negative probability in the framework of combined probability

**Mark Burgin**

University of California, Los Angeles
405 Hilgard Ave.
Los Angeles, CA 90095

**Abstract**

Negative probability has found diverse applications in theoretical physics. Thus, construction of sound and rigorous mathematical foundations for negative probability is important for physics. There are different axiomatizations of conventional probability. So, it is natural that negative probability also has different axiomatic frameworks. In the previous publications (Burgin, 2009; 2010), negative probability was mathematically formalized and rigorously interpreted in the context of extended probability. In this work, axiomatic system that synthesizes conventional probability and negative probability is constructed in the form of combined probability. Both theoretical concepts – extended probability and combined probability – stretch conventional probability to negative values in a mathematically rigorous way. Here we obtain various properties of combined probability. In particular, relations between combined probability, extended probability and conventional probability are explicated. It is demonstrated (Theorems 3.1, 3.3 and 3.4) that extended probability and conventional probability are special cases of combined probability.

## 1. Introduction

All students are taught that probability takes values only in the interval [0,1]. All conventional interpretations of probability support this assumption, while all popular formal descriptions, e.g., axioms for probability, such as Kolmogorov's axioms (Kolmogorov, 1933), canonize this premise. However, scientific research and practical problems brought researchers to the necessity to use more general concepts of probability.

Historically, Hermann Weyl was the first researcher to unconsciously encounter negative probability in (Weyl, 1927). As Groenewold (1946) and Moyal (1949) demonstrated, Weyl considered abstractly and mathematically the transformation that gave what was later called the Wigner distribution, but Weyl did not discuss the fact that it can obtain negative values.

Later Eugene Wigner (1932) introduced a function, which looked like a conventional probability distribution and has later been better known as the Wigner quasi-probability distribution because in contrast to conventional probability distributions, it took negative values, which could not be eliminated or made nonnegative. The importance of Wigner's discovery for foundational problems was not recognized until much later. Another outstanding physicist, Dirac (1942) not only supported Wigner's approach but also introduced the physical concept of negative energy. He wrote:

"*Negative energies and probabilities should not be considered as nonsense. They are well-defined concepts mathematically, like a negative of money*."

Pauli (1956) gave examples of generation of negative probabilities by the renormalization procedure, demonstrating how mathematical models of quantum fields introduced the split of such fields into a positive frequency part and a negative frequency part, giving up thereby the relativistic invariance and also the symmetry of particles and antiparticles.

Dirac (1974) described other useful physical interpretations of negative probabilities.

Discussing computer simulation of quantum reality, Richard Feynman in his keynote talk on Simulating Physics with Computers (1950) said:

*"The only difference between a probabilistic classical world and the equations of the quantum world is that somehow or other it appears as if the probabilities would have to go negative … "*

Later Feyman (1987) wrote a special paper on negative probability where he discussed different examples demonstrating how negative probabilities naturally come to physics and beyond.

Many physicists, including Heisenberg (1931), Wigner (1932), Dirac (1942), Pauli (1943; 1956), and Feynman (1950; 1987), used negative probabilities in their research.

Scully, et al (1994) showed how the micromaser which-path detector eliminated these conceptual difficulties, emphasizing that the concept of negative probability yields useful insight into the Einstein-Podolsky-Rosen (EPR) problem. In essence, negative probabilities appear again and again in different domains of physics.

In such a way, negative probabilities a little by little have become a popular although questionable technique in physics. Many physicists started using negative probabilities to solve various physical problems (cf., for example, (Sokolovski and Connor, 1991; Scully, et al, 1994; Youssef, 1994; 1995; 2001; Han, *et al*, 1996; Curtright and Zachos, 2001; Sokolovski, 2007; Bednorz and Belzig, 2009; Hofmann, 2009)). For instance, Sutherland (2005) examines the usual formalism describing probability density for a particle's position in order to understand the meaning of negative values for probability. He demonstrates how backwards-in-time effects provided a meaningful interpretation for the notion of negative probability. Schneider (2005-2007) demonstrates how Bell's Theorem leads to predictions of negative probabilities. Hartle (2008) developed the quantum mechanics of closed systems, such as the universe, using extended probabilities, which incorporate negative probabilities. Rothman and Sudarshan (2001) demonstrate how quantum mechanics predicts and even necessitates appearance of negative probabilities.

Following Dirac (1974), Finkelstein (2010) writes that negative probabilities have a useful physical interpretation and are no more problematical than a negative bank deposit or a negative energy. In a similar way, the change in population of a system is positive for increase (input) and negative for decrease (output). Finkelstein argues that his approach permits physicists to use a single vector space *V* for both input and output operations, instead of *V* and its dual *V*<sup>*</sup>, making possible to distinguish input from output statvectors in *V* by the signs of their norms.

Moreover, Han, Hwang and Koh (1996) demonstrated the necessity of negative probability measures for some quantum mechanical predictions. Cereceda (2000) extended these results, proving the necessity of negative probability measures for local hidden-variable models in quantum mechanics. Ferrie and and Emerson (2008) also give arguments for necessity of negative probability in quantum theories.

From physics, probabilities came to other disciplines. For instance, they are also used in mathematical models in finance. The concept of *risk-neutral* or *pseudo-probabilities*, i.e., probabilities that can take negative values, is a popular concept in finance and has been numerously applied, for example, in credit modeling by Jarrow and Turnbull (1995), and Duffie and Singleton (1999). Haug (2004; 2007) extended the risk-neutral framework to allow negative probabilities, demonstrating how negative probabilities can help add flexibility to financial modeling. Burgin and Meissner (2010; 2012) applied negative probabilities to mathematical finance developing more flexible models, which adequately represent not only normal functioning of financial systems but also situations of a financial crisis, such as the 2007-2008 worldwide financial crisis or 2011-2013 European financial crisis.

Negative probabilities have been also used in machine learning (Lowe, 2004/2007). Based on mathematical ideas from classical statistics and modern ideas from information theory, Lowe (2004/2007) also argues that the use of non-positive probabilities is both inevitable and natural.

Mathematical problems of negative probabilities have been also studied by Bartlett, Allen, Khrennikov and Burgin. Bartlett (1945) worked out the mathematical and logical consistency of negative probabilities demonstrating how to accurately operate with them but without establishing rigorous foundation. Allen (1976) made an important contribution to understanding negative probabilities suggesting an axiom system for *signed probability theory*, which included negative probabilities. However, as Allen wrote himself, his system was essentially incomplete and was not sufficiently rigorous, for example, in defining random events although his axioms were correct and consistent. Important mathematical results for negative probabilities were obtained by Khrennikov, who writes in his book published in 2009, that it provides the first mathematical theory of negative probabilities. However his theory does not include the classical probability theory because it is developed not in the conventional setting of real numbers pivotal for physics or economics but in the framework of *p*-adic analysis. The mathematical grounding for the negative probability in the real number domain was developed by Burgin (2009) who constructed a Kolmogorov type axiom system building a mathematical theory of extended probability as a probability function, which is defined for random events and can take both positive and negative real values. As a result, extended probabilities include negative probabilities. It was also demonstrated that the classical probability is a positive section (fragment) of extended probability.

There are also meticulous overviews of various advanced concepts of probability. One was written by Mückenheim in 1980s (Mückenheim, 1986)). Another recently written by Burgin reflects the contemporary situation in this area (Burgin, 2012).

In this paper, we employ the axiomatic approach to negative probability. Although Kolmogorov's axioms for probability (1933) remain the most popular in probability theory, other axiom systems for probability were also constructed by Bernstein (1917), Ramsey (1931), de Finetti (1937), Cox (1946), Savage (1954), Watson and Buede (1987) and some other researchers.

Thus, it is natural that negative probability also has different axiomatic frameworks. One of such systems was constructed and studied by Burgin (2009; 2010) in the context of extended probability. In this work, axiomatic system that synthesizes conventional probability and negative probability defines combined probability. Both theoretical concepts – extended probability and combined probability – stretch conventional probability to negative values in a mathematically rigorous way. Various properties of combined probability are obtained in Section 2. Relations of the new approach to extended probability and conventional probability are explicated in Section 3. Interpretations of combined probability are considered in Section 4.

Here we study the symmetric case - symmetric combined probability. The general case of combined probability is studied elsewhere.

## 1.1 Constructions and notation

We remind that if $X$ is a set, then $|X|$ is the number of elements in (cardinality of) $X$ (Kuratowski and Mostowski, 1967). If $A \subseteq X$, then the complement of $A$ in $X$ is defined as $C_X A = X \setminus A$. When $X$ is always fixed, we write simply $X \setminus A = CA$. The set of all subsets of a set $X$ is denoted by $2^X$.

A system **B** of sets is called a *set ring* (Kolmogorov and Fomin, 1989) if it satisfies conditions (R1) and (R2):

(R1)  $A, B \in \mathbf{B}$ implies $A \cap B \in \mathbf{B}$.

(R2)  $A, B \in \mathbf{B}$ implies $A \Delta B \in \mathbf{B}$ where $A \Delta B = (A \setminus B) \cup (B \setminus A)$.

For any set ring **B**, we have $\emptyset \in \mathbf{B}$ and $A, B \in \mathbf{B}$ implies $A \cup B, A \setminus B \in \mathbf{B}$.

Indeed, if $A \in \mathbf{B}$, then by R1, $A \setminus A = \emptyset \in \mathbf{B}$. If $A, B \in \mathbf{B}$, then $A \setminus B = ((A \setminus B) \cup (B \setminus A)) \cap A = (A \Delta B) \cap A \in \mathbf{B}$ by conditions (R1) and (R2).

If $A, B \in \mathbf{B}$ and $A \cap B = \emptyset$, then $A \Delta B = A \cup B \in \mathbf{B}$. Consequently, $A \cup B = (A \setminus B) \cup (B \setminus A) \cup (A \cap B) \in \mathbf{B}$ as $A \cap B \in \mathbf{B}$ by R1 and $(A \setminus B) \cup (B \setminus A) = A \Delta B \in \mathbf{B}$ by R2 because $(A \setminus B) \cap (B \setminus A) = \emptyset$ and $(A \setminus B) \cup (B \setminus A) \cap (A \cap B) = \emptyset$. Thus, a system **B** of sets is a set ring if and only if it is closed with respect union, intersection and set difference.

The set **CI** of all closed intervals $[a, b]$ in the real line ***R*** is a set ring.

The set **OI** of all open intervals $(a, b)$ in the real line ***R*** is a set ring.

A set ring **B** with a unit element, i.e., an element $E$ from **B** such that for any $A$ from **B**, we have $A \cap E = A$, is called a *set algebra* (Kolmogorov and Fomin, 1989).

The set **BCI** of all closed subintervals of the interval $[a, b]$ is a set algebra.

The set **BOI** of all open subintervals of the interval $[a, b]$ is a set algebra.

A set algebra **B** closed with respect to complement is called a *set field*.

## 2. Elements of combined probability theory

To define combined probability, we need some concepts and constructions, which are described below.

Let us consider a set $\Omega$, which consists of two irreducible parts (subsets) $\Omega^+$ and $\Omega^-$, i.e., neither of these parts is equal to its proper subset, and a set **F** of subsets of $\Omega$. Elements from the set $\Omega^+$ are called *provisional elementary positive events*, while subsets of the set $\Omega^-$ are called *provisional elementary negative events* or *provisional elementary antievents*. Subsets of the set $\Omega^+$ are called *provisional positive events*, while subsets of the set $\Omega^-$ are called *provisional negative events* or *antievents*. Subsets of the set $\Omega$ have a common name a *provisional event*. Thus, there are *provisional positive events*, which contain only elements from $\Omega^+$ and there are *provisional negative events*, which contain only elements from $\Omega^-$. All other subsets of the set $\Omega$ are *provisional mixed events*.

There are different examples of negative events (antievents). They are usually connected to negative objects. For instance, encountering a negative object is a negative event.

A transparent example of negative objects is given by antiparticles, which form an antimatter counterpart of quantum particles. For instance, the antiparticle of the electron $e^-$ is called the positron and denoted by $e^+$. When a particle meets its antiparticle, they annihilate each other and disappear, their combined rest energies becoming available to appear in other forms. This kind of annihilation is reflected by a formal operation in the theory of extended probability (Burgin, 2009).

When Paul Dirac extended quantum mechanics to include special relativity, he derived a formula known as the Dirac equation. He noticed in 1928 that this equation predicted that an electron should have a positively charged counterpart (Dirac, 1930). Although the majority of physicists opposed Dirac's conclusion, this hypothetical particle, the positron, was soon discovered in the cosmic radiation by Carl Anderson in 1932.

Another example of negative events (negative objects) is given by antipatterns, which are negative design patterns in the software industry field (Koenig, 1995; Laplante and Neill, 2005). A software design pattern is an *antipattern* if it is a repeated pattern of action, process or structure that initially appears to be beneficial, but ultimately produces more bad consequences than beneficial results, and for which an alternative solution exists that is clearly documented, proven in actual practice and repeatable. Happening of an antipattern is a negative event.

Note that annihilation occurs not only in physics where particles and antiparticles annihilate one another, but also in ordinary life of people. For instance, a person has stocks of two companies. If in 2009, the first set of stocks gave profit $1,000, while the second set

of stocks dropped by $1,000, then the combined income was $0. The loss annihilated the profit.

Here is even a simpler example. A person finds $20 and looses $10. As the result, the amount of this person's money has increased by $10. The loss annihilated part of the gain.

Thus, it is natural to assume that if $A$ belongs to **F**, then it does not contain pairs $\{w, -w\}$ with $w \in \Omega$ and $\mathbf{F}^+ = \{X \in \mathbf{F}; X \subseteq \Omega^+\}$ is a set algebra (cf., for example, (Kolmogorov and Fomin, 1989)) with respect to union with annihilation and $\Omega^+$ is a member of $\mathbf{F}^+$.

Elements from **F**, i.e., subsets of $\Omega$ that belong to **F**, are called *random events*.

Elements from $\mathbf{F}^+ = \{X \in \mathbf{F}; X \subseteq \Omega^+\}$ are called *positive random events*.

Elements from $\Omega^+$ that belong to $\mathbf{F}^+$ are called *elementary positive random events* or simply, *elementary positive random events*.

If $w \in \Omega^+$, then $-w$ is called the *antievent* of $w$. We assume that $-(-w) = w$.

Elements from $\Omega^-$ that belong to $\mathbf{F}^-$ are called *elementary negative random events* or *elementary random antievents*.

Note that the model is symmetric because positive events are, in this sense, antievents of the corresponding negative events, e.g., $E$ is the antievent of $-E$.

For any set $X \subseteq \Omega$, we define
$$X^+ = X \cap \Omega^+,$$
$$X^- = X \cap \Omega^-,$$
$$-X = \{-w; w \in X\}$$

and
$$\mathbf{F}^- = \{-A; A \in \mathbf{F}^+\}$$

If $A \in \mathbf{F}^+$, then $-A$ is called the *antievent* of $A$.

Elements from $\mathbf{F}^-$ are called *negative random events* or *random antievents*.

Now it is possible to give an axiomatic definition of a probability function for combined probability. Here we treat only the symmetric finite case when $\Omega = \{w_1, w_2, w_3, \ldots, w_n, -w_1, -w_2, -w_3, \ldots, -w_n\}$, $\Omega^+ = \{w_1, w_2, w_3, \ldots, w_n\}$ and $\Omega^- = \{-w_1, -w_2, -w_3, \ldots, -w_n\}$. It means that each positive elementary event $w$ has its opposite $-w$ and each negative elementary event $-v$ has its opposite $-(-v) = v$.

**Definition 2.1.** The function $p$ from **F** to the set **R** of real numbers is called a *combined probability function*, if it satisfies the following axioms:

**CP 1** (*Algebraic structure*). **F** is a set algebra that has $\Omega$ as its element.

**CP 2** (*Normalization*). For any event $A \in \mathbf{F}$, $-1 \leq p(A) \leq 1$.

**CP 3** (*Finite additivity*)

$$p(A \cup B) = p(A) + p(B)$$

for all events $A, B \in \mathbf{F}$ such that

$$A \cap B = A \cap -B = \emptyset$$

**CP 4** (*Symmetry of events*). If $A \in \mathbf{F}$, then $-A \in \mathbf{F}$.

**CP 5** (*Symmetry of probability*). $p(A) = -p(-A)$ for any random event $A$ from **F**.

Axiom CP5 allows the following interpretation: the negative probability $-p$ of an event $A$ means that it is possible to expect the opposite event $-A$ with the probability $p$. Thus, it is allows us to understand the negative probability of an event $A$ as the positive probability of the opposite event $-A$. For instance, when the probability of winning $100 is equal to $-\frac{3}{4}$, it means that the probability loosing $100 is equal to $\frac{3}{4}$.

Note that if all elementary events belong to **F**, then Axiom CP5 is equivalent to the following property.

**CP 5a.** $p(w) = -p(-w)$ for any elementary event $A$.

Taking a combined probability function $p$, we define the following sets:

$\Omega_{+p} = \{ w \in \Omega; p(w) > 0 \}$ is the set of all *positively evaluated* elementary events;

$\Omega_{-p} = \{ w \in \Omega; p(w) < 0 \}$ is the set of all *negatively evaluated* elementary events;

$\Omega_{0p} = \{ w \in \Omega; p(w) = 0 \}$ is the set of all *neutrally evaluated* elementary events.

Assuming that all elementary events belong to **F**, we have

$$\Omega = \Omega_{+p} \cup \Omega_{-p} \cup \Omega_{0p}$$

In general, we call a combined probability function $p$ *digitalized* when all elementary events belong to **F**.

It is possible to ask whether the subsets of these sets form a set algebra.

**Proposition 2.1.** The sets $2^{\Omega_{+p}}$, $2^{\Omega_{-p}}$ and $2^{\Omega_{0p}}$ of all subsets of the sets $\Omega_{+p}$, $\Omega_{-p}$ and $\Omega_{0p}$, correspondingly, are set algebras for an arbitrary combined probability function $p$.

Indeed, if $A, B \in \Omega_{+p}$, then $A \cap B$, $A \cup B$, $A \setminus B \in \Omega_{+p}$ and $\Omega_{+p}$ is the unit element with respect intersection $\cap$. Thus, $2^{\Omega_{+p}}$ is a set algebra.

The same is true for the sets $2^{\Omega_{-p}}$ and $2^{\Omega_{0p}}$.

We also have the following sets:

$\mathbf{F}_{+p} = \{ A \in \mathbf{F}; p(A) > 0\}$ is the set of all *positively evaluated* random events;

$\mathbf{F}_{-p} = \{ A \in \mathbf{F}; p(A) < 0\}$ is the set of all *negatively evaluated* random events;

$\mathbf{F}_{0p} = \{ A \in \mathbf{F}; p(A) = 0\}$ is the set of all *neutrally evaluated* random events.

It is possible to ask whether these sets, $\mathbf{F}_{+p}$, $\mathbf{F}_{-p}$ and $\mathbf{F}_{0p}$, are set algebras. As the following examples show, the answer is negative.

**Example 2.1**. Let us consider random events $A = \{w, v, u\}$ and $B = \{a, b, u\}$ with $p(w) = p(v) = p(-u) = p(a) = p(b) = 1/5$. Then $p(A) = 2/5 - 1/5 = 1/5$ and $p(B) = 2/5 - 1/5 = 1/5$. Thus, $A$ and $B$ belong to $\mathbf{F}_{+p}$. At the same time, $A \cap B = \{u\}$ and thus, $p(A \cap B) = p(u) = -1/5$. So, $A \cap B$ does not belong to $\mathbf{F}_{+p}$.

**Example 2.2**. Let us consider random events $A = \{w, v, u\}$ and $B = \{w, v, z\}$ for which $p(w) = p(v) = p(-z) = p(-u) = 1/5$. Then $p(A) = 2/5 - 1/5 = 1/5$ and $p(B) = 2/5 - 1/5 = 1/5$. Thus, $A$ and $B$ belong to $\mathbf{F}_{+p}$. At the same time, $A \cup B = \{ w, v, u, z \}$ and thus, $p(A \cup B) = 2/5 - 2/5 = 0$. So, $A \cup B$ does not belong to $\mathbf{F}_{+p}$. Besides, $A \setminus B = \{u\}$ and thus, $p(A \setminus B) = p(u) = -1/5$. So, $A \setminus B$ also does not belong to $\mathbf{F}_{+p}$.

Given a random event $A$, we define the following sets:

$A_{+p} = A \cap \Omega_{+p}$ is the *positively valuated* part of $A$;

$A_{-p} = A \cap \Omega_{-p}$ is the *negatively valuated* part of $A$;

$A_{0p} = A \cap \Omega_{0p}$ is the *neutrally valuated* part of $A$.

For instance, in the case of extended probability (Burgin, 2009), we have $\Omega_{+p} \subseteq \Omega^+$ and $\Omega_{-p} \subseteq \Omega^-$. If the probability of any elementary event is larger than zero, then we have equalities instead of inclusions of sets.

Axiom CP5 implies the following results.

**Proposition 2.2.** An event $A$ belongs to $\mathbf{F}$ if and only if the event $-A$ belongs to $\mathbf{F}$.

**Lemma 2.1**. For any random event $A$, if $p(A) = 0$, then $p(-A) = 0$ and if $p(A) \neq 0$, then $p(-A) \neq 0$.

**Corollary 2.1**. If $p(w) = 0$, then $p(-w) = 0$ and if $p(w) \neq 0$, then $p(-w) \neq 0$.

Iterative application of Axiom CP3 gives us the following result.

**Lemma 2.2.** If $A_i \cap A_j = A_i \cap -A_j = \emptyset$ for all $i \neq j$, $i, j = 1, 2, 3, \ldots, n$, then

$$p(\bigcup_{i=1}^{n} A_i) = \sum_{i=1}^{n} p(A_i)$$

Proof is performed by induction and application of Axiom CP3.

Axiom CP5 also implies the following result.

**Proposition 2.3.** If $A \in \mathbf{F}$ and $A = -A$, then $p(A) = 0$.

Indeed, if $p(A) = a$, then by Axiom CP5, $p(-A) = -a$. However, by the initial conditions, $A = -A$. Consequently, $a = -a$. This is possible only when $a = 0$.

As $\emptyset = -\emptyset$, we have the following result.

**Corollary 2.2.** $p(\emptyset) = 0$.

As $\Omega = -\Omega$, we have the following result.

**Corollary 2.3.** $p(\Omega) = 0$.

Note that in a non-symmetric case, $p(\Omega)$ is not necessarily equal to zero although $p(\emptyset)$ is always equal to zero.

Let us show that the condition from Axiom CP3 is symmetric.

**Lemma 2.3.** If $-A \cap B = \emptyset$, then $A \cap -B = \emptyset$.

*Proof.* Let us assume that $-A \cap B = \emptyset$ but $A \cap -B \neq \emptyset$. Then there is an elementary event $w$ such that $w \in A$ and $w \in -B$. Thus, $-w \in -A$ and $-w \in B$.

This contradicts the condition $-A \cap B = \emptyset$ and completes the proof.

Definitions imply the following results.

**Lemma 2.4.** For any events $A$ and $B$, we have $-(A \cap B) = -A \cap -B$.

**Lemma 2.5.** For any elementary event $w$, we have $w \in \Omega_{+p}$ if and only if $-w \in \Omega_{-p}$.

**Corollary 2.4.** $\Omega_{-p} = -\Omega_{+p}$ and $p(\Omega_{-p}) = -p(\Omega_{+p})$.

The setting of extended probability implies annihilation of opposite events (Burgin, 2009). The premises of combined probability do not demand annihilation of opposite events but such events do not influence the final outcome and may be excluded.

**Definition 2.2.** a) An event $A$ is called *reduced* if it does not contain subsets of the form $\{w, -w\}$.

b) The event *A* obtained by exclusion all subsets of the form {*w*, -*w*} event *A* is called the *reduction* of *A* and is denoted by R*A*.

c) The event $D = A \setminus RA$ is called the *reducible part* of *A* and is denoted by C*A*.

Let us consider events *A* and *B*.

**Lemma 2.6.** If *B* is reduced and $A \subseteq B$, then *A* is also reduced.

**Corollary 2.5.** If *A* is reduced, then $A \cap B$ and $A \setminus B$ are also reduced for any event *B*.

**Remark 2.1.** Even when both *A* and *B* are reduced, their union is not always reduced.

Let us consider an event *A*.

**Lemma 2.7.** The following conditions are equivalent:

1) The event *A* is reduced.

2) $A \cap -A = \emptyset$.

3) $A^+ \cap -A^- = \emptyset$.

4) $A_{+p} \cap A_{-p} = \emptyset$ and $A_{0p}$ is reduced.

Indeed, a pair {*w*, -*w*} belongs to *A* if and only if it also belongs to -*A*. So, if $A \cap -A = \emptyset$, then there no such pairs in *A* and it is reduced. Besides, if *A* is reduced, then it does not contain pairs of the form {*w*, -*w*}, and thus, $A \cap -A = \emptyset$. This gives us property (1). Properties (2) – (4) are proved in a similar way.

**Lemma 2.8.** For any event *A*, $RA = A \setminus (A \cap -A)$.

As $-(-A) = A$, Lemma 2.8 implies the following results.

**Corollary 2.6.** For any event *A*, $R(-A) = (-A) \setminus (A \cap -A)$.

**Corollary 2.7.** An event *A* is reduced if and only if -*A* is reduced.

There are connection between reductions of an event *A* and the event -*A*.

**Lemma 2.9.** For any event *A*, $RA = -(R(-A))$.

Reduction is not a monotone operation, i.e., it is possible that for some events *A* and *B*, $A \subseteq B$ but $RA \not\subseteq RB$. Indeed, taking $A = \{u, w, -w\}$ and $B = \{u, w\}$, we see that $A \subseteq B$ but $RA = \{u\} \not\subseteq RB = \{u, w\}$.

**Proposition 2.4.** If a combined probability function is digitalized, then for any random event *A* from **F**, we have

$$p(A) = \sum_{w \in A} p(w) \qquad (1)$$

*Proof.* For any event $A$, we have $A = RA \cup CA$. As $RA \cap CA = RA \cap -CA = \emptyset$, by Axiom CP3, we have $p(A) = p(RA) + p(CA)$. As $p(w) = -p(-w)$ for any elementary event $w$,

$$p(CA) = \sum_{w \in CA} p(w) = 0 \qquad (2)$$

Applying Axiom CP3 and induction on the number of elementary events in $RA$, we obtain

$$p(RA) = \sum_{w \in RA} p(w) = 0 \qquad (3)$$

Combining formulas (2) and (3), we obtain formula (1).

Proposition is proved.

**Proposition 2.5.** If a combined probability function $p$ is digitalized, then $\mathbf{F} = 2^{\Omega}$.

Indeed, any set algebra is closed with respect to union and any element from $2^{\Omega}$ is the union of elementary events.

**Corollary 2.8.** If $A \in \Omega_{+p}$, then $p(A) > 0$.

**Corollary 2.9.** If $A \in \Omega_{-p}$, then $p(A) < 0$.

**Corollary 2.10.** $p(A_{+p}) > 0$.

**Corollary 2.11.** $p(A_{-p}) < 0$.

Operation of reduction implies equivalence between events.

**Definition 2.3.** Two events $A$ and $B$ are *equivalent*, $A \approx B$, if they have the same reduction, i.e., $RA = RB$.

**Proposition 2.6.** If $A \approx B$, then $p(A) = p(B)$ for all sets $A, B \in \mathbf{F}$.

*Proof.* Let us assume that the event $A$ contains a pair of elementary events $w$ and $-w$. Then taking the event $D = A \setminus \{w, -w\}$, we have

$$D \cap \{w, -w\} = \emptyset$$

and

$$D \cap -\{w, -w\} = D \cap \{w, -w\} = \emptyset$$

Thus, by Axiom CP3, we have

$$p(A) = p(D) + p(\{w, -w\})$$

At the same time, by Axiom CP5,

$$p(\{w, -w\}) = -p(-\{w, -w\}) = -p(\{w, -w\})$$

Consequently, $p(\{w, -w\}) = 0$ and $p(A) = p(D)$.

In such a way, we can exclude all opposite pairs $\{w, -w\}$ from both events $A$ and $B$ without changing their combined probability. As they are equivalent, the event obtained after all exclusions will be the same for both of them. Thus, $p(A) = p(B)$.

Proposition is proved.

This result means that we do not need to demand annihilation of opposite events as it is done for extended probability (Burgin, 2009) because opposite events in $A$ do not influence the combined probability $p(A)$ of any random event $A$.

As $A \approx RA$, we have the following result.

**Corollary 2.12.** For any random event $A$, $p(A) = p(RA)$.

**Lemma 2.10.** If $A_{0p}$ and all its subsets belong to $\mathbf{F}$, then $p(A_{0p}) = 0$.

Indeed, in this case by Axiom CP3, $p(A_{0p}) = \sum_{w \in A_{0p}} p(w) = 0$ because $p(w) = 0$ for all $w \in A_{0p}$.

**Lemma 2.11.** For any event $A$, $R(CA) = R(-A)$.

*Proof.* If $w \in -A$, then either $-w \in -A$ or $-w \notin -A$. In the first case, the pair $\{w, -w\}$ also belongs to $A$ but is eliminated from $-A$ by reduction. Thus, neither $w$ nor $-w$ belongs to $CA$ and as a matter of fact, to $R(CA)$. It means that all elementary events that are reduced from $-A$ do not belong to $R(CA)$.

In the second case, $-w \in A$ and thus, $w$ belongs to $CA$ and cannot be eliminated from $CA$. Consequently, all elementary events that belong to $R(-A)$ are not eliminated from $CA$ in reduction. Besides, all elementary events that do not belong either to $A$ or to $-A$ form pairs $\{u, -u\}$, which are eliminated from $CA$ by reduction. Consequently, $R(CA) = R(-A)$.

Lemma is proved.

**Corollary 2.13.** For any random event $A$, $p(CA) = p(-A)$.

This result shows that even when we do not know the whole space $\Omega$, it is possible to find the combined probability of the complement of a random event.

Now let us consider operations with events. In general, the union of two reduced events is not a reduced event. For instance, taking reduced events $A = \{w\}$ and $B = \{-w\}$, we see that their union $A \cup B = \{w, -w\}$ is not a reduced event. Thus, it is useful to use the *reduced union* of random events, which is defined as

$$A \cup_R B = R(A \cup B)$$

Many properties of reduced union are similar to properties of conventional union. For instance, we have the following result.

**Proposition 2.7.** For any events *A*, *B* and *C*, the following equalities are valid:

(a) $A \cup_R B = B \cup_R A$

(b) $A \cup_R (B \cup_R C) = (A \cup_R B) \cup_R C$

However, some properties of reduced union and conventional union are different. For instance, conventional union is a monotone operation, i.e., if $C \subseteq B$, then $A \cup C \subseteq A \cup B$ for any event (set) *A*. For reduced union, this is not true in general. For instance, taking reduced events $A = \{w, v, u\}$, $B = \{-w, -v\}$ and $C = \{-w\}$, we see that $C \subseteq B$, $A \cup_R B = \{u\}$, while $A \cup_R C = \{u, v\}$ and thus, $A \cup_R C$ is larger than $A \cup_R B$.

Another new operation with events is *double difference*, which is defined in the following way:

$$A \setminus\!\setminus B = A \setminus (B \cup -B)$$

Lemma 2.5 implies the following result.

**Corollary 2.14.** If *A* is reduced, then $A\setminus\!\setminus B$ is also reduced for any event *B*.

There are certain connections between probabilities of reduced random events.

**Proposition 2.8.** If random events *A* and *B* are reduced, then

$$p(A \cup B) = p(A) + p(B) - p(A \cap B).$$

*Proof.* Let us take two reduced random events *A* and *B*. By properties of sets, we have

$$A = (A \cap B) \cup (A \cap -B) \cup (A \setminus\!\setminus B)$$

and

$$B = (B \cap A) \cup (B \cap -A) \cup (B \setminus\!\setminus A)$$

At first, we see that $(A \setminus\!\setminus B) \cap (A \cap B) = (A \setminus\!\setminus B) \cap (A \cap -B) = \varnothing$ because $A \setminus\!\setminus B = A \setminus (B \cup -B)$. Besides, $(A \setminus\!\setminus B) \cap -(A \cap B) = (A \setminus\!\setminus B) \cap -(A \cap -B) = \varnothing$ because the event *A* is reduced. Indeed, if $w \in (A \setminus\!\setminus B) \cap -(A \cap B)$, then $w \in (A \setminus\!\setminus B)$ and $w \in -(A \cap B)$. Consequently, $w \in A$ and $w \in -A$. Thus, $-w \in A$ and this contradicts our assumption that the event *A* is reduced. By the same argument, $(A \setminus\!\setminus B) \cap -(A \cap -B) = \varnothing$.

In addition, $(A \cap B) \cap (A \cap -B) = (A \cap B) \cap -(A \cap -B) = \varnothing$. Indeed, if $w \in (A \cap B) \cap (A \cap -B)$, then $w \in (A \cap B)$ and $w \in (A \cap -B)$. Consequently, $w \in B$ and $w \in -B$. Thus, $-w \in B$ and this contradicts our assumption that the event *B* is reduced. In the same way, if *w*

∈ (A ∩ B) ∩ -(A ∩ -B), then w ∈ (A ∩ B) and w ∈ -A ∩ B. Consequently, w ∈ A and w ∈ -A. Thus, -w ∈ A and this contradicts our assumption that the event A is reduced.

Then it is possible to apply Lemma 2.2, getting

$$p(A) = p(A \cap B) + p(A \cap -B) + p(A \backslash\backslash B) \qquad (4)$$

In a similar way, we obtain

$$p(B) = p(B \cap A) + p(B \cap -A) + p(B \backslash\backslash A) \qquad (5)$$

As $A \cap B = B \cap A$ and $p(A \cap -B) = p(A \cap -B)$ by Axiom CP5, we have

$$p(A) + p(B) = 2p(A \cap B) + p(A \backslash\backslash B) + p(B \backslash\backslash A) \qquad (6)$$

At the same time,

$$A \cup B = (A \cap B) \cup (B \setminus A) \cup (A \backslash\backslash B) =$$
$$= (A \cap B) \cup (B \backslash\backslash A) \cup (B \cap -A) \cup (A \backslash\backslash B) \cup (A \cap -B)$$

Besides, $(B \backslash\backslash A) \cap (A \backslash\backslash B) = (A \cap B) \cap (B \backslash\backslash A) = (A \cap B) \cap (A \backslash\backslash B) = \emptyset$ by the properties of sets. We have also proved that

$$(A \cap B) \cap (A \cap -B) = (A \cap B) \cap (-A \cap B) =$$
$$(A \backslash\backslash B) \cap (-A \cap B) = (A \backslash\backslash B) \cap (A \cap -B) = \emptyset$$

By a similar argument, we have

$$(B \backslash\backslash A) \cap (-A \cap B) = (B \backslash\backslash A) \cap (A \cap -B) = \emptyset.$$

In addition, $(A \cap -B) \cap (-A \cap B) = \emptyset$. Indeed, if $w \in (-A \cap B) \cap (A \cap -B)$, then $w \in (-A \cap B)$ and $w \in (A \cap -B)$. Consequently, $w \in B$ and $w \in -B$. Thus, $-w \in B$ and this contradicts our assumption that the event B is reduced. Then it is possible to apply Lemma 2.2, getting

$$p(A \cup B) = p(A \cap B) + p(A \backslash\backslash B) + p(B \backslash\backslash A) + p(B \cap -A) + p(A \cap -B)$$

However, $B \cap -A = -(A \cap -B)$. Thus, by Axiom CP5, we have $p(B \cap -A) = -p(A \cap -B)$ and

$$p(A \cup B) = p(A \cap B) + p(A \backslash\backslash B) + p(B \backslash\backslash A) \qquad (7)$$

Comparing equalities (6) and (7), we come to the necessary result

$$p(A \cup B) = p(A) + p(B) - p(A \cap B) \qquad (8)$$

Proposition is proved.

Formula (8) is the same for conventional probability but for combined probability the proof is more complicated.

For a conventional probability $P$, if for random events $A$ and $B$, $A \subseteq B$ holds, then $p(A) \leq p(B)$. This is not true for combined probabilities in a general case. Instead, we have a more sophisticated property. Let us assume that a combined probability function $p$ digitalized, i.e., all elementary events from $\Omega$ belong to $\mathbf{F}$.

**Proposition 2.9.** If $A_{+p} \subseteq B_{+p}$ and $B_{-p} \subseteq A_{-p}$, then $p(A) \leq p(B)$.

Indeed, there are more elementary random events with positive probability evaluation in $B$ than in $A$ and there are less elementary random events with negative probability evaluation in $B$ than in $A$.

**Definition 2.4.** a) A random event $A$ is *positively valuated* by $p$ if $p(A) > 0$.
b) A random event $A$ is *negatively valuated* by $p$ if $p(A) < 0$.

For instance, the event $A_{+p}$ is positively valuated by $p$ and the event $A_{-p}$ is negatively valuated by $p$ for any random event $A$.

**Proposition 2.10.** All random events from $2^{\Omega_{+p}}$ are positively valuated by $p$ and the combined probability function $p$ is monotone on these events.

*Proof.* The first part of the statement is established in Corollary 2.9. So we need to prove only the second part. Because all elementary events from $2^{\Omega_{+p}}$ are evaluated by $p$, we have $2^{\Omega_{+p}} \subseteq \mathbf{F}$.

Let us take random events $A$ and $C$ from $2^{\Omega_{+p}}$ such that $C \subseteq A$. As $2^{\Omega_{+p}}$ is a set algebra, the set $D = A \setminus C$ belongs to $2^{\Omega_{+p}}$ and by construction $D \cap C = \emptyset$. Assume that the random event $A$ is reduced, then $D \cap \text{-}C \subseteq A \cap \text{-}A = \emptyset$. Thus, $D \cap \text{-}C = \emptyset$ and by Axiom CP3, $p(A) = p(C) + p(D)$. Consequently, $p(A) \geq p(C)$ since $p(D) \geq 0$. Note that equality holds only when $C = A$.

When the random event $A$ is not reduced, we take events $RA$ and $RC$ and show that $p(RA) \geq p(RC)$. As $p(A) = p(RA)$ and $p(C) = p(RC)$, we obtain $p(A) \geq p(C)$. Note that it is possible that $RC = RA$ even when $C \neq A$.

Proposition is proved.

**Proposition 2.11.** All random events from $2^{\Omega_{-p}}$ are negatively valuated by $p$ and the combined probability function $p$ is antitone on these events.

Indeed, the first part is established in Corollary 2.10, while the proof of the second part is similar to the proof of Proposition 2.8.

There are different kinds of combined probability functions.

**Definition 2.5.** a) A combined probability function $p$ is called *positively normalized* if there is a random event $A$ such that $p(A) = 1$.

b) A combined probability function $p$ is called *negatively normalized* if there is a random event $D$ such that $p(D) = -1$.

c) A combined probability function $p$ is called *normalized* if it is both positively normalized and negatively normalized.

In contrast to extended probability (Burgin, 2009; 2012), we do not include the normalization condition in the list of basic axioms for combined probability because there are situations in real life when we do not know all possible outcomes of an experiment and thus, we do not have even a single random event probability of which is equal to 1. Examples of such situations are considered, for example, in (Montagna, 2012).

In a similar way, completeness is not included in the list of basic axioms for combined probability.

**Definition 2.6.** a) The space $\Omega$ is *positively complete* for $p$ if $p(\Omega_{+p}) = 1$.

b) The space $\Omega$ is *negatively complete* for $p$ if $p(\Omega_{-p}) = -1$.

There are definite relations between normalization and completeness.

Axiom CP5 implies the following result.

**Proposition 2.12.** a) If for any random event $A$, $p(A) > 0$ implies that any elementary event $w$ from $A$ is positively evaluated by $p$, then a combined probability function $p$ is positively normalized if and only if $\Omega$ positively complete for $p$.

b) If for any random event $A$, $p(A) < 0$ implies that any elementary event $w$ from $A$ is negatively evaluated by $p$, then a combined probability function $p$ is negatively normalized if and only if $\Omega$ negatively complete for $p$.

*Proof.* By definition if $p(\Omega_{+p}) = 1$, then $p$ is positively normalized.

At the same time, let us take a positively normalized combined probability function $p$. Then there is a random event $A$ such that $p(A) = 1$. It is possible to assume that the event $A$ is reduced.

As all elementary events from $\Omega_{+p}$ are positively evaluated by $p$, they all belong to **F**. Consequently (cf. Proposition 2.4), $\Omega_{+p} \in$ **F**. The event $\Omega_{+p}$ is reduced because if $p(w) > 0$, then $p(-w) < 0$.

By the initial condition, $A \subseteq \Omega_{+p}$ and by Lemma 2.10, $A \subseteq R\Omega_{+p}$ because $A$ is reduced. Thus, $\Omega_{+p} = A \cup D$ where $D = \Omega_{+p} \setminus A$ and $p(\Omega_{+p}) = p(A) + p(D)$. As all elementary events from $D$ are positively evaluated, $p(D) > 0$. However, by Axiom CP2, $p(\Omega_{+p}) \leq 1$. Thus, $\Omega_{+p} = A$ and $p(\Omega_{+p}) = 1$.

Part (a) is proved.

The proof of part (b) is similar.

Proposition is proved.

The initial conditions in Proposition 2.10 are essential as the following example demonstrates.

**Example 2.3**. Let us consider the space of events $\Omega = \{w, v, u, -w, -v, -u\}$ such that $p(w) = -p(-w) = 1/3$, $p(\{w, v, u\}) = -p(\{-w, -v, -u\}) = 1$. Then the combined probability function $p$ is positively normalized although $\Omega_{+p} = \{w\}$ and $p(\Omega_{+p}) = 1/3$. It is also negatively normalized although $\Omega_{-p} = \{-w\}$ and $p(\Omega_{-p}) = -1/3$.

**Proposition 2.13.** A combined probability function $p$ is positively normalized if and only if it is negatively normalized.

Indeed, if $\Omega_{+p} \in$ **F**, then $\Omega_{-p} \in$ **F**. Besides, $\Omega_{+p} = -\Omega_{-p}$ and by Axiom CP5, $p(\Omega_{+p}) = -p(\Omega_{-p})$. So, $p(\Omega_{+p}) = 1$ if and only if $p(\Omega_{-p}) = -1$.

**Corollary 2.15.** A combined probability function $p$ is positively normalized if and only if it is normalized.

**Remark 2.2.** Proposition 2.6 shows that in the symmetric case studied here, the concepts of a positively normalized combined probability function, of a negatively normalized combined probability function and of a normalized combined probability function coincide. However, for a non-symmetric probability space, this is not always true.

We show that Axiom EP of extended probability (Burgin, 2009) can be proved for combined probability of reduced random events.

**Proposition 2.14** (*Decomposition property*). For any reduced random event $A \in$ **F**, if $A_{+p}, A_{-p}, A_{0p} \in$ **F**, then

$$p(A) = p(A_{+p}) + p(A_{-p})$$

*Proof.* By construction, $A = A_{+p} \cup A_{-p} \cup A_{0p}$. Definitions imply that $A_{0p} \cap A_{-p} = A_{0p} \cap A_{+p} = A_{+p} \cap A_{-p} = \emptyset$ and $A_{0p} \cap -A_{-p} = A_{0p} \cap -A_{+p} = \emptyset$ because by Lemma 2.1, if $p(w) = 0$, then $p(-w) = 0$ and if $p(w) \neq 0$, then $p(-w) \neq 0$. Besides, $A_{+p} \cap -A_{-p} = \emptyset$. Indeed, if $w \in A_{+p}$ and $w \in -A_{-p}$, then $-w \in A_{-p}$. This implies that $w, -w \in A$. As $A$ is a reduced random event, it is impossible.

Thus, we can apply Lemma 2.2, which gives us

$$p(A) = p(A_{+p}) + p(A_{-p}) + p(A_{0p})$$

As by Lemma 2.6, $p(A_{0p}) = 0$, we have

$$p(A) = p(A_{+p}) + p(A_{-p})$$

Proposition is proved.

Note that for arbitrary events the Decomposition property is not true.

## 3. Combined probability versus conventional probability and extended probability

There are also other ways for building a probability measure that takes negative values. Another approach to this problem has the form of combined probability, which includes negative probability. In this approach, we do not need symmetry of the space $\Omega$ of elementary events although in some cases, it may be symmetric, i.e., containing the negative event $-w$ for each elementary event $w$ and containing the positive elementary event $-u$ for each negative elementary event $u$.

For extended probability, the statement of Axiom CP5 is deduced from the axioms (Burgin, 2009). We cannot do the same for combined probability because the combined probability of a positive event can be negative, i.e., Axiom EP8 is not valid for combined probability. On the other hand, Axiom EP7 can be deduced from Axiom CP5 and Axiom CP3 (cf. Proposition 2.9).

There are inherent relations between extended probability and combined probability.

**Theorem 3.1.** Any extended probability is a combined probability in which all random events are reduced.

*Proof.* As in extended probability theory opposite events annihilate one another, we need only to show that any extended probability function satisfies Axioms CP1 – CP5, or more exactly, that Axioms EP1 – EP8 of extended probability (Burgin, 2009) imply Axioms CP1 – CP5 of combined probability.

At first, we recall axioms of extended probability.

**EP 1** (*Order structure*). There is a graded involution $\alpha: \Omega \to \Omega$, i.e., a mapping such that $\alpha^2$ is an identity mapping on $\Omega$ with the following properties: $\alpha(w) = -w$ for any element $w$ from $\Omega$, $\alpha(\Omega^+) \supseteq \Omega^-$, and if $w \in \Omega^+$, then $\alpha(w) \notin \Omega^+$.

**EP 2** (*Algebraic structure*). $\mathbf{F}^+ \equiv \{X \in \mathbf{F}; X \subseteq \Omega^+\}$ is a set algebra that has $\Omega^+$ as a member.

**EP 3** (*Normalization*). $P(\Omega^+) = 1$.

**EP 4** (*Composition*) $\mathbf{F} \equiv \{X; X^+ \in \mathbf{F}^+ \ \& \ X^- \in \mathbf{F}^- \ \& \ X^+ \cap -X^- \equiv \emptyset \ \& \ X^- \cap -X^+ \equiv \emptyset\}$.

**EP 5** (*Finite additivity*)

$$P(A \cup B) = P(A) + P(B)$$

for all sets $A, B \in \mathbf{F}$ such that

$$A \cap B \equiv \emptyset$$

**EP 6** (*Annihilation*). $\{v_i, w, -w; v_i, w \in \Omega \ \& \ i \in I\} = \{v_i; v_i \in \Omega \ \& \ i \in I\}$ for any element $w$ from $\Omega$.

**EP 7**. (*Adequacy*) $A = B$ implies $P(A) = P(B)$ for all sets $A, B \in \mathbf{F}$.

**EP 8**. (*Non-negativity*) $P(A) \geq 0$, for all $A \in \mathbf{F}^+$.

Now we can use these axioms to deduce Axioms CP1 – CP5.

Axiom CP1 is proved in Theorem 2 from (Burgin, 2009).

Axiom CP2 is proved in Proposition 11 from (Burgin, 2009).

Axiom CP3: As by Proposition 8 from (Burgin, 2009), $P(A) = -P(-A)$ for any random event $A$ from $\mathbf{F}$, Axiom CP3 is directly implied by Axiom EP5.

Axiom CP4 is deduced from Axiom EP4.

Axiom CP5 is proved in Proposition 8 from (Burgin, 2009).

Theorem is proved.

This result allows researchers to infer many properties of extended probability from properties of combined probability.

The converse of Theorem 3.1 is not true. For instance, let us take $\Omega = \{w, -w\}$ and define $p(w) = -1$ and $p(-w) = 1$. Then $p$ is a combined probability function but it is not and extended probability function because in it a positive event has the negative probability.

However, in some cases, combined probability can coincide with extended probability.

**Theorem 3.2.** If a combined probability function $p$ is normalized, all random events are reduced, $\Omega_{+p} \cup \Omega_{0p} \supseteq \Omega^+$ and $\Omega_{-p} \cup \Omega_{0p} \supseteq \Omega^-$, then $p$ is an extended probability function.

*Proof.* We need to show that when the initial conditions are satisfied, the combined probability function $p$ satisfies Axioms EP1 – EP8, or more exactly, that Axioms CP1 – CP5 of combined probability imply Axioms EP1 – EP8 of extended probability.

Axiom EP1 is true because we consider the symmetric case of combined probability taking $\mathbf{F} = 2^{\Omega}$.

Axiom EP2 is true when we take $\mathbf{F}^+ = 2^{\Omega^+}$.

Axiom EP3 is true because the combined probability function $p$ is normalized and by the initial conditions, all subsets from $\Omega^+$ are non-negatively valuated by $p$.

Axiom EP4 is true because $\mathbf{F}^+ = 2^{\Omega^+}$, $\mathbf{F}^- = 2^{\Omega^-}$ and $A = A^+ \cup A^-$ for any $A \subseteq \Omega$.

Axiom EP5: Let us consider events $A, B \in \mathbf{F}$ such that $A \cap B \equiv \emptyset$. By the initial conditions, the probability function $p$ is digitalized. Thus, by Proposition 2.3, we have

$$p(A) = \sum_{w \in A} p(w)$$

$$p(B) = \sum_{w \in B} p(w)$$

and

$$p(A \cup B) = \sum_{w \in A \cup B} p(w) = \sum_{w \in A} p(w) + \sum_{w \in B} p(w)$$

because $A$ and $B$ do not have common elements. Thus, $P(A \cup B) = P(A) + P(B)$, i.e., Axiom EP5 is true.

Axiom EP6 is true because opposite elementary events do not change probability (Proposition 2.5).

Axiom EP7 follows from Proposition 2.5.

Axiom EP8 follows from the condition $\Omega_{+p} \cup \Omega_{0p} \supseteq \Omega^+$ and the equality $\mathbf{F}^+ = 2^{\Omega^+}$.

Theorem is proved.

This result shows that under definite conditions, combined probability becomes extended probability

**Theorem 3.3.** Conventional probability is a completely asymmetric, i.e., without negative and negatively valuated events, positively normalized combined probability.

*Proof.* We need to show that when there are no negatively valuated events, a positively normalized combined probability function $p$ satisfies axioms K1 – K3 of Kolmogorov for conventional probability.

At first, we remind Kolmogorov's axioms for probability (Kolmogorov, 1933):

**K 1**. (*Non-negativity*) $P(A) \geq 0$, for all $A \in \mathbf{F}$.

**K 2**. (*Normalization*) $P(\Omega) = 1$.

**K 3**. (*Finite additivity*)

$$P(A \cup B) = P(A) + P(B)$$

for all sets $A, B \in \mathbf{F}$ such that

$$A \cap B = \emptyset.$$

Now we can use Axioms CP1 – CP5 to deduce Axioms K1 – K3.

Axiom K1 is true because all events are positively valuated by $p$.

Axiom K2 is true because $p$ is positively normalized.

Axiom K3 follows from Axiom CP3 because opposite events do not exist in this case.

Theorem is proved.

This result shows that combined probability is a natural extension of conventional probability and it is possible to obtain many properties of conventional probability as corollaries of properties of combined probability..

**Theorem 3.4.** The restriction $p_+$ of a combined probability function $p$ on the set $\mathbf{F}_{+p}$ is a conventional probability function, i.e., it satisfies Kolmogorov axioms, if and only if $p$ is positively normalized and $\mathbf{F}_{+p}$ is a set algebra.

*Proof. Necessity.* If the restriction $p_+$ of a positively normalized combined probability function $p$ on the set $\mathbf{F}_{+p}$ satisfies Kolmogorov axioms, then by Axiom K2, $p$ is positively normalized, while $\mathbf{F}_{+p}$ is a set algebra by the initial conditions for the convenient probability function (Kolmogorov, 1933), the set of all random events, in our case $\mathbf{F}_{+p}$, is a set algebra.

*Sufficiency*. We need to show that when $\mathbf{F}_{+p}$ is a set algebra, then a positively normalized combined probability function *p* satisfies axioms K1 − K3 of Kolmogorov for conventional probability.

Axiom K1 is true because all events from $\mathbf{F}_{+p}$ are positively valuated by *p*.

Axiom K2 is true because *p* is positively normalized.

Axiom K3 follows from Axiom CP3 because opposite events do not belong to $\mathbf{F}_{+p}$.

Theorem is proved.

**Remark 3.1.** Both conditions from Theorem 3.4 are essential. Indeed, Examples 2.1 and 2.2 show that $\mathbf{F}_{+p}$ is not always a set algebra. Besides, there exist combined probability functions *p* that are not positively normalized. For instance, taking an arbitrary combined probability functions *p*, we can define a new combined probability functions *q* by the formula $q(A) = ½ p(A)$ for any random event *A* from $\mathbf{F}$. We can see that all axioms CP1 − CP5 remain true for *q*. At the same time, even when *p* is positively normalized, *q* is not positively normalized.

**Conclusion**

Extended probability studied in (Burgin, 2009; 2010; Burgin and Meissner, 2010; 2012) describes randomness as if we live in a perfect world where positive events have positive probabilities, while negative events have negative probabilities. Combined probability studied here reflects randomness in a more realistic way assuming that any event can have either positive probability or negative probability depending on external conditions.

Obtained results show that combined probability is a natural extension of both conventional probability and extended probability. This allows researchers to obtain many properties of conventional probability and extended probability as corollaries of properties of combined probability.

Here only the first step of combined probability studies is presented. It is necessary to explore combined probability for infinite and continuous probability spaces $\Omega$, to study combined probability distributions and to investigate moments of random variables, such as expectation and variance, with respect to combined probability.

It would be also interesting to find other axiomatic systems for combined probability, which are similar to the axiomatic systems for conventional probability constructed by Ramsey (1931), de Finetti (1937), Cox (1946) and Savage (1954).